\documentclass[11pt]{amsart}

\theoremstyle{definition}

\numberwithin{equation}{section}

\begin{document}

\title
{Fast proof of functional equation for $\zeta(s)$}

\author{Luis B\'{a}ez-Duarte}

\date{13 May 2003}
\maketitle

The functional equation for $\zeta(s)$ seems certain to remain a most important
topic. More than a few dozen proofs for it are already known; perhaps this is a lower
bound only. In \textit{R\'{e}coltes et Semailles\textit} Grothendieck admonishes
mathematicians to look anew at old concepts in solitude and in absolute,
childlike innocence. Following his dictum I found a very short proof of the functional
equation, only to realize later -one expects to be forgiven- that it was an even
shorter version of Titchmarsh's (\cite{titch} p. 15) first proof of the functional
equation based on the fundamental Fourier transform view of
$\zeta(s)$ that appears, easily overlooked, on page 14 formula
(2.1.5) of his treatise. Namely, with $\rho(x)$ standing for the fractional part of
$x$,

\begin{equation}\label{mellin}
\frac{\zeta(s)}{-s}=\int_0^\infty x^{-s-1}\rho(x)dx, \ \ \ (0<\Re s <1),
\end{equation}
which, by the way, is the basis for the Nyman-Beurling approach \cite{beurling},
\cite{baez} to the Riemann hypothesis.

Starting from (\ref{mellin}) we can write
\begin{eqnarray}\nonumber
(2^s-1)\frac{\zeta(s)}{s}&=&\int_0^\infty x^{-s-1}(\rho(x)-\rho(2x))dx\\\nonumber
&=&\int_0^\infty x^{-s-1}
\sum_{n=1}^\infty \frac{1}{n\pi}(\sin 4n\pi x-\sin 2n\pi x)dx.
\end{eqnarray}
valid for $0<\Re s<1$.
On account of the uniform boundedness of the partial sums we can interchange integral
and sum to obtain

\begin{eqnarray}\nonumber
(2^s-1)\frac{\zeta(s)}{s}&=&
\sum_{n=1}^\infty \frac{1}{n\pi}\int_0^\infty x^{-s-1}(\sin 4n\pi x-\sin 2n\pi
x)dx\\\nonumber
&=&-(2^s-1)2^s \pi^{s-1}\Gamma(-s)\sin\frac{\pi s}{2}\zeta(1-s),
\end{eqnarray}
where the last equation is valid at first when $-1<\Re s<0$, but it obviously provides
the analytic extension to the whole plane.

\bibliographystyle{amsplain}

\begin{thebibliography}{12}
\bibitem{baez}
L. B\'{a}ez-Duarte, \textit{A strengthening of the Nyman-Beurling 
criterion for the Riemann hypothesis}, preprint 2002, posted in arxiv
math.NT/0205003, to appear in Atti Accad. Naz. Lincei Cl. Sci. Fis. Mat. Natur.
Rend. Lincei (9) Mat. Appl. 1.

\bibitem{beurling}
A. Beurling, \textit{A closure problem related to the Riemann Zeta-function}, Proc. Nat. Acad.
Sci. \textbf{41} (1955), 312-314.
\bibitem{nyman}
B. Nyman, \textit{On some groups and semigroups of translations}, Thesis, Uppsala,
1950.
\bibitem{titch}
E. C. Titchmarsh, \textit{The Theory of the Riemann-Zeta Function}, Clarendon Press,
Oxford, 1951.
\end{thebibliography}

\ \\
\noindent Luis B\'{a}ez-Duarte\\
Departamento de Matem\'{a}ticas\\
Instituto Venezolano de Investigaciones Cient\'{\i}ficas\\
Apartado 21827, Caracas 1020-A, Venezuela
\ \\
\email{lbaezd@cantv.net}

\end{document}